\documentclass[letterpaper, 10 pt, conference]{ieeeconf}  % Comment this line out
\usepackage{amsmath} % assumes amsmath package installed
\usepackage{amssymb}  % assumes amsmath package installed

\IEEEoverridecommandlockouts                              % This command is only
\usepackage{color}
\usepackage{multirow}
\usepackage{graphicx}
\usepackage{subfig}
\usepackage{epstopdf}
\newtheorem{definition}{Definition}
\newtheorem{theorem}{Theorem}
\newtheorem{lemma}{Lemma}
\newtheorem{remark}{Remark}
\newtheorem{assumption}{Assumption}
\newtheorem{proposition}{Proposition}

\newcommand{\argmin}{\mathop{\mathrm{argmin}}\limits} 

\newcommand{\Id}{ \mathrm{Id}}
\newcommand{\bR} { {\mathbb R}}
\newcommand{\bN} { {\mathbb N}}
\newcommand{\mH} { {\mathcal H}}
\newcommand{\fix} { {\mathrm{fix}}}

\title{\LARGE \bf
Towards time-varying proximal dynamics \\ in Multi-Agent Network Games 
}

\author{Carlo Cenedese$^{1}$ \and Yu Kawano$^{1}$ \and Sergio Grammatico$^{2}$ \and Ming Cao$^{1}$% <-this % stops a space
\thanks{$^{1}$ Engineering and Technology Institute Groningen (ENTEG), Faculty of Science and Engineering, University of Groningen, The Netherlands         {\tt\small c.cenedese@rug.nl} , {\tt\small y.kawano@rug.nl} and {\tt\small m.cao@rug.nl}. The work of Cenedese, Kawano and Cao was supported in part by the European Research Council (ERC-CoG-771687) and the Netherlands Organization for Scientific Research (NWO-vidi-14134).
}%
\thanks{$^{2}$ Delft Center for Systems and Control, TU Delft, The Netherlands
{\tt\small s.grammatico@tudelft.nl}. The work of Grammatico was partially supported by NWO under research projects OMEGA (613.001.702) and P2P-TALES (647.003.003).
}}

\begin{document}

\maketitle
\thispagestyle{empty}
\pagestyle{empty}

% ABSTRACT %%%%%%%%%%%%%%%%%%%%%%%%%%%%%%%%%%%%%%%%%%%%%%%%%%%%%%%%%%%%%%%%%%%%%%%%%%%%%%%
\begin{abstract}
Distributed decision making in multi-agent networks has recently attracted significant research attention thanks to its wide applicability, e.g.\ in the management and optimization of computer networks, power systems, robotic teams, sensor networks and consumer markets. Distributed decision-making problems can be modeled as inter-dependent optimization problems, i.e., multi-agent game-equilibrium seeking problems, where noncooperative agents seek an equilibrium by communicating over a network. To achieve a network equilibrium, the agents may decide to update their decision variables via proximal dynamics, driven by the decision variables of the neighboring agents. In this paper, we provide an operator-theoretic characterization of convergence with a time-invariant communication network. For the time-varying case, we consider adjacency matrices that may switch subject to a dwell time. We illustrate our investigations using a distributed robotic exploration example.
\end{abstract}

% INTRODUCTION %%%%%%%%%%%%%%%%%%%%%%%%%%%%%%%%%%%%%%%%%%%%%%%%%%%%%%%%%%%%%%%%%%%%%%%%%%%%%%%
\section{Introduction}

\subsection{Motivation: Multi-agent decision making over networks}
Multi-agent decision making over networks is currently a vibrant research area in the systems-and-control community, with application in several relevant domains,  such as smart grids \cite{dorfer:simpson-porco:bullo:16, grammatico:18tcns}, traffic and information networks \cite{jaina:walrand:10, barrera:garcia:15}, social networks \cite{Ghaderi2014, etesami:basar:15}, consensus and flocking groups \cite{olfati-saber:murray:04, olfati-saber:06}, robotic and sensor networks \cite{martinez:bullo:cortes:frazzoli:07, sankovic:johansson:stipanovic:12}.

In distributed computation and communication, the main advantage is that each decision maker, in short, \textit{agent}, can keep its own data private and exchange information with selected agents only.  
Essentially, in networked multi-agent systems, the state (or decision) variables of each agent evolve as a result of \textit{local decision making}, e.g. local constrained optimization, and \textit{distributed communication} with some neighboring agents, via a communication graph. Typically, the aim of the agents is reaching a collective equilibrium state, where no agent can benefit from further updating its state variables.

\subsection{Literature overview: Multi-agent optimization and multi-agent network games}

Multi-agent dynamics for solving a set of inter-dependent optimization problems arise naturally from distributed optimization and distributed equilibrium seeking in network games.
Multi-agent convex constrained optimization has been widely studied in the last decade: in \cite{nedic:ozdaglar:parrillo:10} with uniformly bounded subgradients, and either homogeneous constraint sets or time-invariant, complete communication graphs with uniform weights; 
in \cite{lee:nedic:13} with differentiable cost functions with Lipschitz continuous and uniformly bounded gradients; and, more generally, in \cite{falsone:margellos:garatti:prandini:17}, where convergence is proven via \textit{vanishing} step sizes. Network games among agents with convex compact local constraints have been considered before: in \cite{parise:gentile:grammatico:lygeros:15} with strongly convex quadratic cost functions and time-invariant communication graph; in \cite{koshal:nedic:shanbhag:16} \cite{salehisadaghiani:pavel:16}, with differentiable cost functions with Lipschitz continuous gradient, strictly convex cost functions, and undirected, possibly time-varying, communication graphs; and in \cite{grammatico:18tcns} with general local convex cost and quadratic proximal term, time-invariant and time-varying communication graphs, subject to technical restrictions. 
%Multi-agent games with convex compact local and also \textit{coupling} constraints have been considered in \cite{yin:shanbhag:mehta:11} under the assumption of strongly convex twice-differentiable cost functions with bounded gradients, with strictly increasing congestion cost term.
 The common feature in multi-agent optimization and games over networks is the presence of a structured, possibly time-varying, communication graph. Therefore, it is interesting to design multi-agent dynamics that involve distributed computation and structured information exchange.

\subsection{Contribution of the paper}
In this paper, we consider proximal dynamics in multi-agent network games with both time-invariant and time-varying communication graphs. In the time-invariant case, we show that global convergence of proximal dynamics holds if the adjacency matrix of the communication graph, assumed strongly connected, is row stochastic and with strictly-positive diagonal elements. Technically, we extend the convergence result in \cite[Th. 1]{grammatico:18tcns}. The use of a row stochastic matrix is highly relevant in applications: it allows an agent to communicate with its neighbors without requiring an adjustment of the rest of the network. In the time-varying case, we consider switching adjacency matrices  subject to a certain dwell time and we show global convergence of the proximal dynamics under switching with sufficiently large dwell time. 
For testing the derived sufficient conditions, we provide linear matrix inequalities.

\subsection{Organization of the paper}

The paper is organized as follows: Section \ref{sec:motivating_application} presents an illustrative multi-robot exploration scenario; Section \ref{sec:problem_formulation} formalizes the problem setup.
We introduce the convergence result for time-invariant multi-agent proximal dynamics Section \ref{section:static_dynamics} and for time-varying, dwell-time switched, proximal dynamics in Section \ref{section:switching_commuting}. A numerical simulations of the considered dynamics is presented in Section \ref{section:numerical_simulations}. Finally, we conclude the paper in Section \ref{section:conclusions}, where we discuss future research directions. %The proofs are given in Appendix, together with basic definitions and results from operator theory.

\subsection{Basic notation}
The set of real, positive, and non-negative are denoted by $\mathbb{R}$, $\mathbb{R}_{>0}$, $\mathbb{R}_{\geq 0}$, respectively; $\overline{\mathbb{R}}:=\mathbb{R}\cup \{\infty\}$. The set of natural numbers is denoted by $\mathbb{N}$, and for $i,j\in\mathbb{N}$, $i\leq j$, we define $\mathbb{N}[i,j]:=[i,j]\cap \mathbb{N}$. For a square matrix $A \in \bR^{n\times n}$, its transpose is denoted by $A^\top$. The $n \times n$ identity matrix is denoted by~$I_n$. For $x_1,\cdots,x_N\in\mathbb{R}^n$, a collective vector $\boldsymbol{x}:=[x_1^\top,\cdots ,x_N^\top ]^\top\in\mathbb{R}^{nN}$ is simply described as $\boldsymbol{x}=[x^1;\cdots;x^N]$.  For two matrices~$A$ and $B$, $A\otimes B$ denotes their Kronecker product. For vectors $x,y \in \bR^n$ and a symmetric and positive definite $n\times n$ matrix $ Q \succ 0$, the weighted inner product and norm are denoted by $\langle x | y\rangle_{Q}$ and $\lVert x \rVert_{Q}$, respectively; the induced matrix norm is denoted by $\lVert A\rVert_{Q}$. For $Q=I_{n}$, the standard inner product, Euclidean norm, and Frobenius norm are obtained. A real $n$ dimensional Hilbert space obtained by endowing $\mathcal H=(\bR^n,\lVert\cdot\rVert)$  with the product $\langle x | y\rangle_{Q}$ is denoted by $\mathcal{H}_Q$.

\subsection{Operator-theoretic notation}
For a function $f:\bR^n\rightarrow\overline{\mathbb{R}}$, define $\mathrm{dom}(f):=\{x\in\bR^n|f(x)<+\infty\}$. 
The subdifferential $\partial f:\mathrm{dom}(f)\rightrightarrows\mathcal{H}_Q$ is defined by $\partial f(x)=\{ u\in \mathcal{H}_Q | (\forall y\in\mathcal{H}_Q)\: \langle y-x|u\rangle +f(x)\leq f(y) \}$.
The proximal operator $\mathrm{prox}_f^Q(x):\mathcal{H}_Q\rightarrow\mathrm{dom}(f)$ is defined by $\mathrm{prox}_f^Q(x):=\mathrm{argmin}_{y\in\bR^n}f(y)+\textstyle\frac{1}{2}\lVert x-y\rVert^2_Q  $.  The resolvent of an operator $A:\mH_Q\rightarrow 2^{\mH_Q}$ is $\mathcal{J}_A :=(\Id+A)^{-1}$.The indicator function $\iota_\mathcal{C}:\bR^n\rightarrow[0,+\infty]$ of $\mathcal{C}\subset \bR$ is defined as $\iota_\mathcal{C}(x)=0$ if $x\in\mathcal{C}$; otherwise $+\infty$. The identity operator is defined by~$\Id$. The Euclidean distance and the  Euclidean distance weighted by $Q$ of a point $x$ to $\mathcal{C}$ are respectively  $\mathrm{d}(x,\mathcal C)$, and $\mathrm{d}_{Q}(x,\mathcal C)$.

%% PROBLEM STATEMENT %%%%%%%%%%%%%%%%%%%%%%%%%%%%%%%%%%%%%%%%%%%%%%%%%

%% PROBLEM FORMULATION %%%%%%%%%%%%%%%%%%%

\section{Motivating, illustrative scenario: Multi-robot exploration}
\label{sec:motivating_application}

\begin{figure}
\centering
\subfloat[][]
{\includegraphics[width=.2\textwidth]{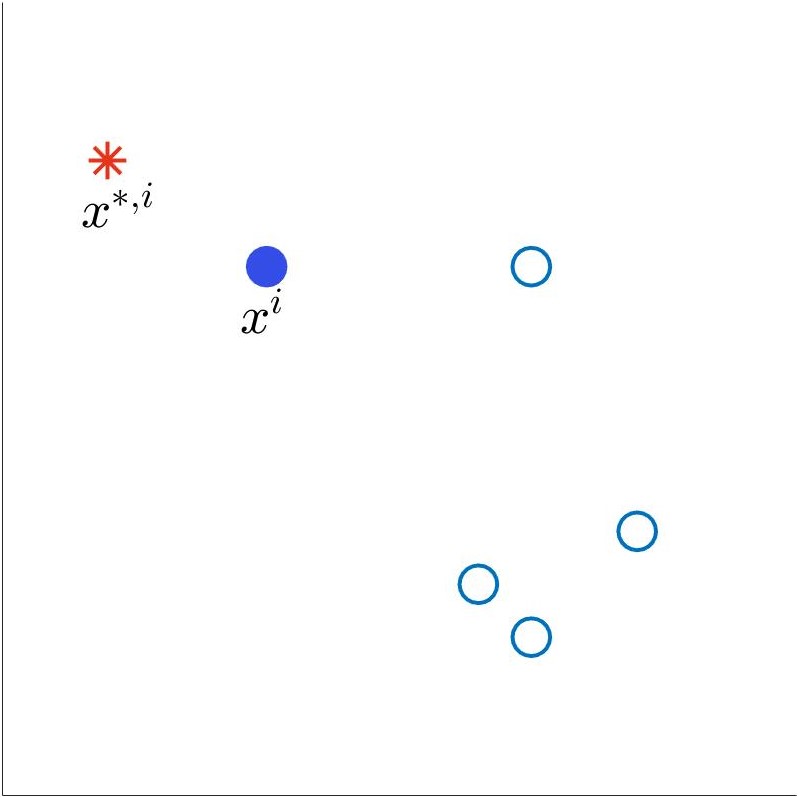}\label{fig:1_a}} \quad %[width=.225\textwidth] % How it was
\subfloat[][]
{\includegraphics[width=.2\textwidth]{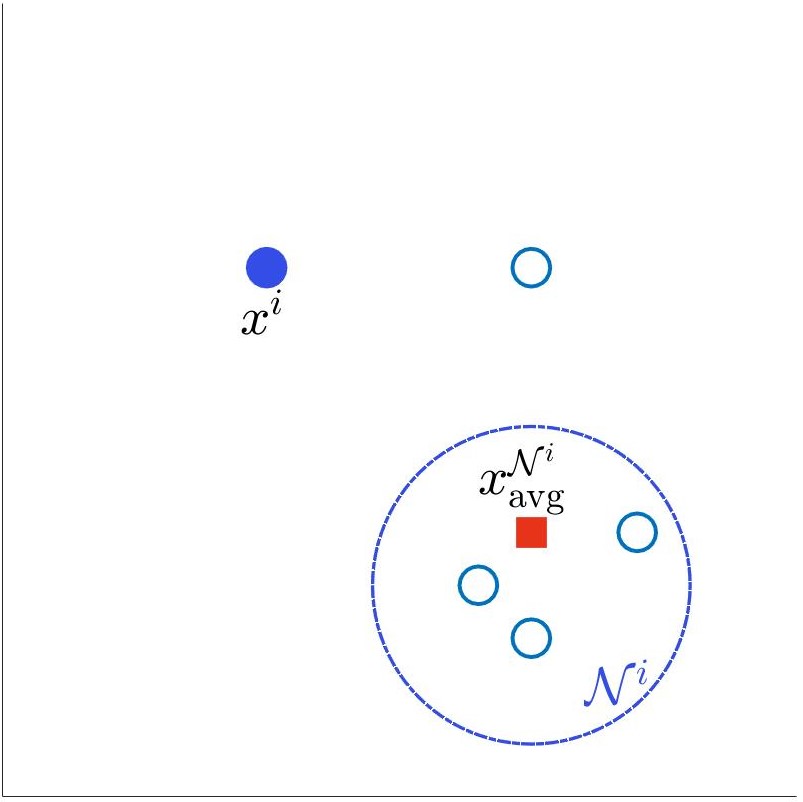}\label{fig:1_b}} \\
\subfloat[][]
{\includegraphics[width=.2\textwidth]{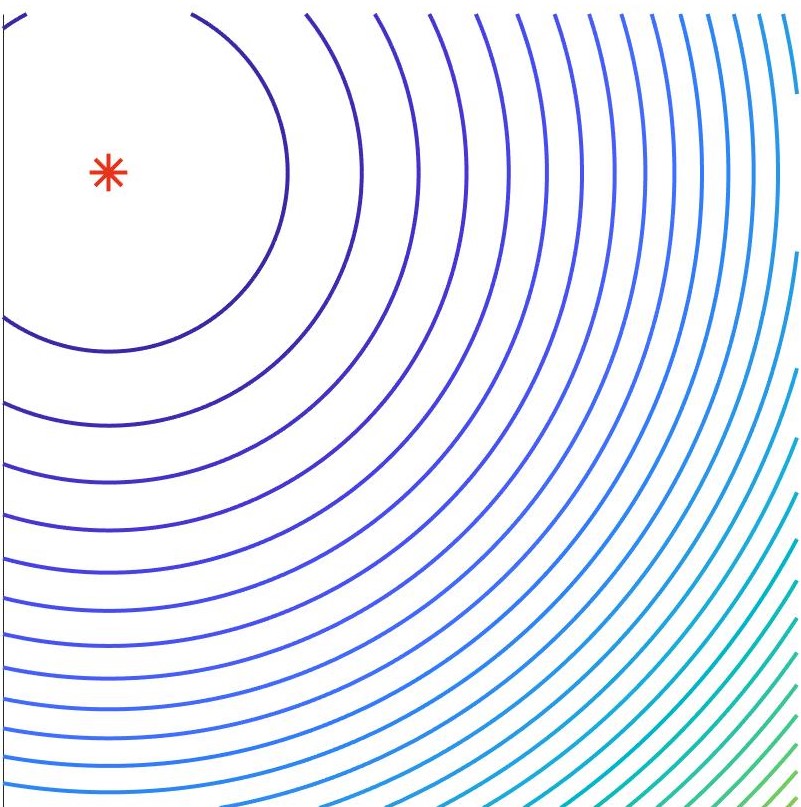}\label{fig:1_c}} \quad
\subfloat[][]
{\includegraphics[width=.2\textwidth]{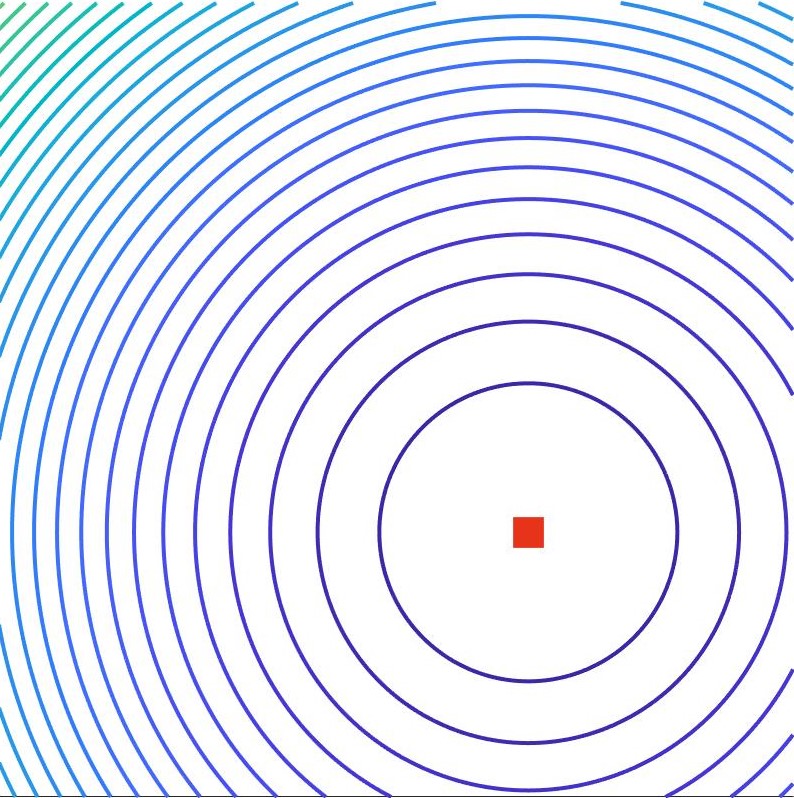}
\label{fig:1_d}}
\caption{(a) Agent $i$ (blue disk) and its target position $x^{*,i}$ (red star), and other robots (light blue circles); (b) $N$ robots, the neighbors $\mathcal{N}^i$ of robot $i$ and $ x^{\mathcal{N}_i}_{\mathrm{avg}}$ (red square), (c) Level sets of the cost function $f^i_{\mathrm{t}}(x^i(k))$; (d) Level sets of the function $g^i( x^i(k),  x^{\mathcal{N}^i}_{\mathrm{avg}})$.}
\label{fig:explenatory_fig_sec2}
\end{figure}

We motivate the paper by the problem of distributed exploration performed by $N$ mobile robots. Let the two-dimensional position of each robot $i\in\bN[1,N]$ at time  $k\in\bN$ be $x^i(k)$ and denote its neighbors indexed by $\mathcal{N}^i$. To each robot $i$, we associate a local cost function that is composed by two separate terms: the local target function $f^i_{\mathrm{t}}(x^i(k))$ and the aggregation term $g^i( x^i(k), x^{\mathcal{N}^i}_{\mathrm{avg}})$, where $ x^{\mathcal{N}^i}_{\mathrm{avg}}$ denotes the weighted averaged positions of its neighboring $\mathcal{N}^i$. The first term penalizes the distance of the robot from its target position $x^{*,i}$, and, by construction, has its minimum at the target position, see Fig.~\ref{fig:1_a},~\ref{fig:1_c}. Instead, the second term penalizes the distance of the position $x^i(k)$ from  $ x^{\mathcal{N}_i}_{\mathrm{avg}}$, hence it plays the role to induce the robots to stay together during their motion, see Fig.~\ref{fig:1_b},~\ref{fig:1_d}. Each robot is assumed to be rational, namely, willing to determine its motion with the aim to minimize its cost function. Overall, the robots shall reach a collective equilibrium state, which we call \emph{network equilibrium}. 

The resulting motion-planning problem can be intended as a game between all the robots involved in the exploration. In fact, the equilibrium points correspond to the trade-offs between the target positions and closeness among robots. For simplicity, in this illustrative example, the collision avoidance between robots is not taken into account. One possible simple structure for the (discrete-time) dynamics of each robot $i$ in the above setup is then
\begin{equation}
x^i(k+1) = \underset{ y\in\mathcal{X}^i}{\mathrm{argmin}}f^i_{\mathrm{t}}(y) + g^i( y, \boldsymbol x^{\mathcal{N}_i})\: ,
\label{eq:general_dynamics_robotic_exploration}
\end{equation}
where the set $\mathcal{X}^i$  represents the motion constraints of the robot.
Whether or not the dynamics in \eqref{eq:general_dynamics_robotic_exploration} will converge to an equilibrium is unclear a-priori, espacially if the set of neighbors, $\mathcal{N}^i(k)$, is time-varying. With this distributed robotic setup in mind, in the following, we  address the convergence problem via an operator-theoretic perspective. 

\section{Technical setup and problem formulation}
\label{sec:problem_formulation}

We consider a network of $N$ agents, where the state of each agent $i\in\mathbb{N} [1,N]$ is denoted by $x^i\in\mathcal{X}^i\subseteq \mathbb{R}^n$ and the set $\mathcal{X}^i$ coincides with the feasible states of agent $i$.
To compute its next state variable, each agent $i$ relies on the states of some neighboring agents. In turn, a network structure arises, described by a weighted digraph. Let us represent the communication between agents by the following $N\times N$ weighted adjacency matrix:
\begin{equation}
P:=[a_{i,j}]=\begin{bmatrix}
a_{1,1} & \cdots & a_{1,N}\\
\vdots & \ddots & \vdots\\
a_{N,1} & \cdots & a_{N,N}
\end{bmatrix},
\label{eq:def_communication_matrix_P}
\end{equation}
where $a_{i,j}\in[0,1]$ is the weight that agent $i$ assigns to the state of agent $j$. If $a_{i,j}=0$, then the state of agent $i$ is independent from that of agent $j$. Furthermore, we assume that each agent $i$ aims at minimizing a cost function $J^i$.

Throughout the paper, we assume compactness and convexity of the local constraint set~$\mathcal{X}^i$ and convexity (not necessarily strict convexity) of the local cost function~$J^i$.    

\smallskip
\begin{assumption}[Local constraints]
\label{ass:local_constrain}
For each agent $i\in\mathbb{N}[1,N]$, the set $\mathcal X^i\subseteq \mathbb{R}^n$ is non-empty, compact and convex.

\hfill\QEDopen
\end{assumption}
\smallskip

\smallskip
\begin{assumption}[Local cost functions]
\label{ass:cost_function}
For each agent $i\in\mathbb{N}[1,N]$, the local cost function $J^i:\mathbb{R}^n\times \mathbb{R}^n\rightarrow\overline{\mathbb{R}}$ is defined by 
\begin{equation}
J^i(y,z) := f_\circ^i(y) + \iota_{\mathcal{X}_i}(y)+\textstyle{\frac{1}{2}}\lVert y-z\rVert^2_{Q_i},
\label{eq:def_costFunction_Ji}
\end{equation}
for some matrix $Q_i\succ0$, where $f^i:=f_\circ^i+\iota_{\mathcal{X}_i}:\mathbb{R}^n\rightarrow\overline{\mathbb{R}}$ is a lower semi-continuous and convex function. 
\hfill\QEDopen
\end{assumption}
\smallskip

In \eqref{eq:def_costFunction_Ji}, the function $f_\circ^i$  is local to agent $i$. For example, it can represent the distance from a desired state. The quadratic term $\textstyle{\frac{1}{2}}\lVert y-z\rVert^2_{Q_i}$ penalizes the distance between the state of agent $i$ and the weighted average among the states of its neighbors. We emphasize that Assumption~\ref{ass:cost_function} requires neither the differentiability of the local cost function, nor the Lipschitz continuity or boundedness of its gradient.

From the above problem setup, since we assume the agents are rational, i.e., willing to minimize their individual cost functions, we consider the following notion of collective equilibrium state, called network equilibrium.

\smallskip
\begin{definition}[Network equilibrium]\label{def:NetworkEquilibrium}
A collective vector $\overline{\boldsymbol x}= [\overline{x}^1;\cdots;\overline{x}^N]\in\mathbb{R}^{nN}$ is a network equilibrium (NWE) if $\forall i\in\mathbb{N}[1,N]$,
\begin{equation}
\overline{x}^i\in  \argmin_{y\in\mathcal{X}^i}J^i(y,\textstyle\sum_{j=1}^Na_{i,j}\overline{x}^j) \:.
\label{eq:NetworkEquilibrium}
\end{equation}
\hfill\QEDopen
\end{definition}
\smallskip

We recall that if there are no self-loops in the adjacency matrix, i.e., $a_{i,i}=0$ for all $i \in\mathbb{N}[1,N]$, then an NWE corresponds to a Nash equilibrium \cite[Remark 1]{grammatico:18tcns}. Under Assumptions~\ref{ass:local_constrain} and~\ref{ass:cost_function}, an NWE always exists. 

The problem studied in this paper is then seeking an NWE (Definition \ref{def:NetworkEquilibrium}), namely, convergence to an NWE from any initial condition. Clearly, in the time-varying case, the definition of NWE is more involved - let us postpone it to Section \ref{section:switching_commuting}.

%% STATIC DYNAMICS %%%%%%%%%%%%%%%%%%%%%%%%%%%%%%%%%%%%%%%%%%%%%%%%%%%%%%%%%%%%%
\section{Time-invariant proximal dynamics}
\label{section:static_dynamics}

As mentioned in the previous section, we assume that each agent is rational and noncooperative. 
Therefore, it is natural to consider the following proximal dynamics for each agent $i\in\mathbb{N}[1,N]$:
\begin{equation}\label{eq:Banach_dynamics_i}
x^i(k+1)  = \mathrm{prox}^{Q_i}_{f^i}\big(\textstyle\sum_{j=1}^Na_{i,j} \, x^j \big),\:   \ \forall k \in \mathbb{N}.
\end{equation}

In the collective vector form, namely, for the collective vector
$$ \boldsymbol{x}(k) := \left[ 
\begin{matrix}
x^1(k);\cdots;x^N(k)
\end{matrix}
\right]\,,
 $$
the dynamics from \eqref{eq:Banach_dynamics_i} read as 
\begin{equation}\label{eq:group_dynamics}
\boldsymbol x(k+1) = \boldsymbol{\mathrm{prox}}_{\boldsymbol{f}}^{\boldsymbol{Q}}( \boldsymbol{A}\,\boldsymbol x(k) ) \,,
\end{equation}
where $\boldsymbol Q$ stands for the block-diagonal matrix
\begin{equation}
\boldsymbol Q := \mathrm{diag}(Q_1, Q_2, \cdots,Q_N) \,,
\label{eq:collective_matrix_Q}
\end{equation}
the matrix $\boldsymbol{A} := P \otimes I_n$ represents the interactions among agents, and the mapping 
$\boldsymbol{\mathrm{prox}}_{\boldsymbol{f}}^{\boldsymbol{Q}}$ is a block-diagonal proximal operator, i.e., 
\begin{equation} \label{eq:group_proximity_operator}
\boldsymbol{\mathrm{prox}}_{\boldsymbol{f}}^{\boldsymbol{Q}}(z) := \mathrm{diag}(\mathrm{prox}_{f^1}^{Q_1}(z^1),\cdots,\mathrm{prox}_{f^N}^{Q_N}(z^N)).
\end{equation}

With the introduced notations,  a collective vector $\overline{\boldsymbol{x}}= [\overline{x}^1;\cdots;\overline{x}^N]$ is an NWE if and only if $\overline{\boldsymbol{x}}\in\mathrm{fix}(\boldsymbol{\textrm{prox}}_{\boldsymbol{f}}^{\boldsymbol{Q}}\circ \boldsymbol{A})$, where $\mathrm{fix}(\cdot)$ stands for the set of fixed points of the operator in its argument. Under Assumptions~\ref{ass:local_constrain} and~\ref{ass:cost_function}, $\mathrm{fix}\left(\boldsymbol{\mathrm{prox}}_{\boldsymbol{f}}^{\boldsymbol{Q}}\circ \boldsymbol{A} \right)$ is non-empty \cite[Th.~4.1.5(b)]{smart1980fixedpoint_theory}, and the convergence problem is well posed.

Therefore, from an operator-theoretic perspective, the proximal dynamics in \eqref{eq:group_dynamics} are the so-called Picard--Banach iteration for the mapping $\boldsymbol{\mathrm{prox}}_{\boldsymbol{f}}^{\boldsymbol{Q}}\circ \boldsymbol{A}$~ \cite[Equ.~1.69]{Bauschke2010:ConvexOptimization}.

%In this paper, w
We assume that the adjacency matrix is row-stochastic with self-loops, and marginally stable, as formalized next.

\smallskip
\begin{assumption}[Row-stochasticity and self-loops]
\label{ass:row_stoch}
The communication graph is strongly connected. The matrix $P=[a_{i,j}]$ in \eqref{eq:def_communication_matrix_P} is row-stochastic, i.e., $a_{i,j} \geq 0$ for all $i,j \in \bN[1,N]$, 
and $\sum_{j=1}^Na_{i,j}=1$ for any $i \in \bN[1,N]$. Moreover, $P$ has strictly-positive diagonal elements, i.e. $\min_{i\in\mathbb{N}[1,N]}a_{i,i}=:\underline{a}>0$.
\hfill\QEDopen
\end{assumption}
\smallskip

%%\smallskip
%%\begin{assumption}[Marginally stable]
%%\label{ass:marginally_stable}
%%The matrix $P=[a_{i,j}]$ in (\ref{eq:def_communication_matrix_P}) is marginally stable. Because of Assumption~\ref{ass:self_loops}, this is equivalent to that the eigenvalue $1$ of $P$ is semi-simple. \hfill\QEDopen 
%%\end{assumption}
%%\smallskip

Now, we are ready to introduce the first result of this paper about the convergence of the proximal dynamics in \eqref{eq:group_dynamics}. 

\smallskip
\begin{lemma}[Global convergence]
\label{th:convergence_banach_dynamics}
Suppose that Assumptions~\ref{ass:local_constrain}--\ref{ass:row_stoch} hold. 
There always exists a matrix $\tilde{\boldsymbol{Q}} \succ 0$ such that, for any $\boldsymbol x(0) \in \mathcal{X}$, the sequence $\left( \boldsymbol x(k) \right)_{k=0}^\infty$ generated by \eqref{eq:group_dynamics} with ${\boldsymbol{Q}} = \tilde{\boldsymbol{Q}}$ converges to an NWE.
\hfill\QEDopen
\end{lemma}
\smallskip

\smallskip
\begin{remark} \label{rem:diagonal}
Lemma \ref{th:convergence_banach_dynamics} extends \cite[Th. 1]{grammatico:18tcns}, since $\boldsymbol A$ is only assumed to be row-stochastic. % where $\tilde{\boldsymbol{Q}}$ is chosen as the identity matrix. 
Note that if the matrix $\tilde{\boldsymbol{Q}}$ can be chosen block-diagonal, then the operator $\boldsymbol{\mathrm{prox}}_{\boldsymbol{f}}^{\tilde{\boldsymbol{Q}}}\circ \boldsymbol{A}$ defines fully distributed dynamics.
\hfill\QEDopen
\end{remark}
\smallskip
    
From the practical point of view, the matrix $\tilde{\boldsymbol{Q}}$ in Lemma \ref{th:convergence_banach_dynamics} can be computed as $\tilde{\boldsymbol{Q}} = \tilde{Q} \otimes I_n$, where the matrix $\tilde{Q}$ solves the following LMI:
\begin{equation}\label{eq:LMI}
A^\top \tilde{{Q}}A \preccurlyeq (2 \eta-1) \tilde{{Q}} + ( 1 - \eta )( {A}^\top \tilde{{Q}} + \tilde{{Q}} {A} )
\end{equation}
for some $\eta \in (0,1)$. We have in fact the following result.

\smallskip
\begin{proposition}
\label{prop:eta_avg_LMI_solvable}
Let $\eta \in (0,1)$. If the LMI in \eqref{eq:LMI} holds, then the sequence $\left( \boldsymbol x(k) \right)_{k=0}^\infty$ generated by \eqref{eq:group_dynamics} with ${\boldsymbol{Q}} = \tilde{{Q}} \otimes I_n$, and $\tilde{{Q}}$ solution to \eqref{eq:LMI}, converges to an NWE.
\hfill\QEDopen
\end{proposition}
\smallskip

It follows from Remark \ref{rem:diagonal} and Proposition \ref{prop:eta_avg_LMI_solvable} that, in order to obtain fully distributed dynamics in \eqref{eq:group_dynamics}, one shall solve the LMI in \eqref{eq:LMI} with diagonal matrix $\tilde{Q}$.

%%\begin{proof}
%%Proving that the matrix $\boldsymbol{A}$ is $\eta-$averaged in $\mathcal{H}_{\tilde{\boldsymbol{Q}}}$, coincides to proving that the auxiliary matrix $R$ in (\ref{eq:R_NE}) is nonexpansive in the same space. Furthermore from Definition~\ref{def:NE_operators} we know that $R$ is nonexpansivness in $\mathcal{H}_{\tilde{\boldsymbol{Q}}}$ if and only if 
%%\begin{equation}
%%R^\top \tilde{\boldsymbol{Q}} R \preceq \tilde{\boldsymbol{Q}}
%%\label{eq:LMI_NE_R}
%%\end{equation}
%%holds for $\tilde{\boldsymbol{Q}}\succ0$. Substituting the definition (\ref{eq:R_NE}) of $R$ in (\ref{eq:LMI_NE_R}) and simple manipulations of the the inequality, lead to the LMI in (\ref{eq:LMI}). This concludes the proof. 
%%\end{proof}

%% SWITCHING WITH COMMUTING MATRICES %%%%%%%%%%%%%%%%%%%%%%%%%%%%%%%%%%%%
\section{Towards time-varying proximal dynamics}
\label{section:switching_commuting}

\subsection{Time-varying setup}
In the previous section, we have assumed that the communication network of the agents is the same for all time instances $k\in \bN$. In practical situations, however, not all agents can update their strategies at the same time instances. More generally, the communication network can change from time to time. To address time-varying scenarios, in this subsection, we consider a time-varying communication matrix, i.e., 
\begin{equation}
P(k):=[a_{i,j}(k)]=\begin{bmatrix}
a_{1,1}(k) & \cdots & a_{1,N}(k)\\
\vdots & \ddots & \vdots\\
a_{N,1}(k) & \cdots & a_{N,N}(k)
\end{bmatrix},
\label{eq:weighted_adjacency_time_dep}
\end{equation}
hence the collective adjacency matrix $\boldsymbol A(k):=P(k)\otimes I_n$.

For simplicity, in the remainder of the paper, we assume that the set of available communication networks is finite.

\smallskip
\begin{assumption}[Finite number of adjacency matrices]
\label{ass:finite_set_P}
There exists $M \in \bN$ such that 
$P(k) \in \mathcal{P} := \left\{ P_1, \ldots, P_M\right\}$ for all $k \in \bN$, where each matrix $P_i\in\mathcal{P}$ satisfies Assumption~\ref{ass:row_stoch}. 
\hfill\QEDopen 
\end{assumption}
\smallskip

%%For instance, if agents update their strategies asynchronously, Assumption~\ref{ass:finite_set_P} holds. Note that even if the original communication is represented by a doubly stochastic matrix, communication matrices for asynchronous updates are not always doubly stochastic but are row stochastic. This means that the conventional results in \cite{Nedic2016:Distributed_Algorithms_Aggregative_Games,Grammatico:CDC_OpinionDynamics} are difficult to extend for analysis of the asynchronously updates, or more generally, the time-varying case. 

To describe the corresponding dynamics, we introduce a switching signal $\sigma:\bN\rightarrow\mathbb{N}[1,M]$ that at each time step $k$ selects an adjacency matrix. Thus, in compact form, we have the following switching dynamics:
\begin{equation}
\boldsymbol x(k+1) = \boldsymbol{\mathrm{prox}}_{\boldsymbol{f}}^{\boldsymbol Q_{\sigma(k)}}\big( \boldsymbol A_{\sigma(k)} \, \boldsymbol x(k) \big) \:,
\label{eq:switching_signal_general}
\end{equation}
where ${\boldsymbol Q_{\sigma(k)}} := Q_{\sigma(k)}\otimes I_n$, and 
$$Q_{\sigma(k)} \in \mathcal{Q}:=\{ Q_1, \ldots, Q_M \succ 0 \}$$ 
for all $k \in \bN$. Moreover, since we are interested in distributed dynamics, we assume that the matrices $Q_i$'s are diagonal.

\smallskip
\begin{assumption} \label{ass:diagonal_Qi}
For each $P_i \in \mathcal{P}$, there exists a diagonal matrix $Q_i \in\mathcal{Q}$ that satisfies Lemma~\ref{th:convergence_banach_dynamics}. 
\hfill $\square$
\end{assumption}
\smallskip

As anticipated in the previous section, for the time-varying case, we need to generalize the concept of an NWE, to what we call a \textit{persistent network equilibrium}.
\smallskip
\begin{definition}[Persistent network equilibrium]
\label{def:persistent_NE}
A persistent network equilibrium (PNWE) is a collective vector in the set 
\begin{equation}\label{eq:persistent_NE}
\mathcal{E}:= \textstyle\bigcap_{ i=1}^{M} \mathrm{fix}\bigl(\boldsymbol{\mathrm{prox}}_{\boldsymbol{f}}^{\boldsymbol Q_{i}}\circ \boldsymbol A_{i} \bigr) \:.
\end{equation}
\hfill\QEDopen
\end{definition}
\smallskip
\begin{assumption}[Existence of PNWE]
\label{ass:persistent_NE}
The set of PNWE is assumed to be non-empty, i.e., $\mathcal E \not = \emptyset$.
\hfill\QEDopen
\end{assumption}
\smallskip
We remind that, in general, even if all $P_i$'s are stable, i.e., correspond to averaged mappings, the switching system in \eqref{eq:switching_signal_general} can be unstable for some switching sequences. Therefore, we shall use tools from switching systems to claim convergence to a PNWE. In the next subsection, we focus on the \textit{dwell-time} approach to establish global convergence.

\subsection{Proximal dynamics with dwell time}
\label{subsec:dwell_time}

With the aim to study global convergence of time-varying proximal dynamics, let us introduce the concept of \textit{dwell-time}~\cite{Liberzon1999:Switching_time,Zhikun_She2017:Dwell_time,Cao:dwell_time}. 

\smallskip
\begin{definition}
A natural number $\tau$ is called a dwell time if the switching times $\bar{k}_1, \bar{k}_2, \dots$ satisfy $\bar{k}_{i+1} - \bar{k}_i > \tau$, for all $i \in \bN$. 
\hfill $\square$
\end{definition}
\smallskip

Before we can establish convergence of the proximal dynamics with switching, we impose two technical assumptions, namely that  all operators are \textit{linearly regular}~\cite[Def.~2.1]{Bausche:convergence_regular_operators} and that the number of switchings is infinite. 

\smallskip
\begin{assumption}
\label{ass:linear_regular}
Assumptions~\ref{ass:local_constrain}, \ref{ass:cost_function}, \ref{ass:finite_set_P} hold and for each $i\in\mathbb{N}[1,M]$, the operator $\boldsymbol{\mathrm{prox}}_{\boldsymbol{f}}^{\boldsymbol Q_i}\circ \boldsymbol A_i$ is linearly regular  on $\mathcal{H}_{\boldsymbol Q_i}$.
\hfill\QEDopen
\end{assumption}
\smallskip
\begin{assumption}[Infinite switching]
\label{ass:infinite_occurence_of_P}
For all $i \in \bN[1,M]$, the switching signal is such that $\sigma(k) = i$ infinitely many times as $k \rightarrow \infty$.
\hfill\QEDopen
\end{assumption}
\smallskip

We are now ready to present a global convergence result for the switching proximal dynamics in \eqref{eq:switching_signal_general} to a PNWE, provided that the dwell time is chosen large enough.

\smallskip
\begin{theorem}[Global convergence under dwell-time]
\label{th:dwell_time_theorem}
Let Assumptions~\ref{ass:persistent_NE}--\ref{ass:infinite_occurence_of_P} hold. 
For any initial condition $\boldsymbol x(0)\in\mathcal{X}$, the sequence $(\boldsymbol x(k))_{k=0}^\infty$ generated by \eqref{eq:switching_signal_general} converges to a PNWE if the dwell time $\tau \in \bN$ is chosen large enough.
\hfill\QEDopen
\end{theorem}
\smallskip
\begin{remark}
A lower bound for the dwell time can be obtained in the form 
\begin{equation}\label{eq:claim_tau}
\tau\geq \tau_{ \min } := \mathrm{log}_{\prod_{j=1}^M \phi_j}\bigg( \frac{1}{2^M}\prod_{j=1}^M\frac{\lambda_{\min,j}}{\lambda_{\max,j}}\bigg) \,,
\end{equation} 
where $\lambda_{\min,j}$ and $\lambda_{\max,j}$ are, respectively, the minimum and the maximum eigenvalues of $Q_j$. The parameters $\phi_1, \ldots, \phi_M > 0$ are described in the Appendix. 
\hfill\QEDopen
\end{remark}

%% NUMERICAL EXAMPLE %%%%%%%%%%%%%%%%%%%%%%%%%%%%%%%%%%%%%%%%%%
\section{Numerical simulations}
\label{section:numerical_simulations}
We resume the setting described in Section~\ref{sec:motivating_application}, namely, the problem of a distributed exploration performed by a network of mobile robots. In the following, we verify the results of Section~\ref{section:static_dynamics} to solve this task in the time-invariant  case.    

%\subsection{Time-invariant dynamics}
%\label{subsec:static_dynamics_obstacle_avoidance}
\subsubsection{Simulation setup}
\label{subsubsec:simulation_setup}
We have $N=4$ agents in the game, where each agent $i\in\bN[1,N]$ is a moving mobile robot, and the state $x^i(k)$ is its position in the plane at time $k\in\bN$, hence $n=2$. The robots are able to move in all directions from their current positions, but the maximum range of movement is limited inside a square $\mathcal{X}^i(k)$, centred in $x^i(k)$ and of edge $r$. The weighted adjacency matrices in \eqref{eq:def_communication_matrix_P} and its collective counterpart are 
$$P=\textstyle \left[\begin{smallmatrix}
\frac{1}{2}& \frac{1}{2}& 0&  0\\
\frac{4}{10}& \frac{1}{2}& \frac{1}{10} &  0\\
\frac{1}{4}&   \frac{1}{4}&  \frac{1}{4}& \frac{1}{4}\\
\frac{1}{4}&   \frac{1}{4}&  \frac{1}{4}& \frac{1}{4}  
\end{smallmatrix}\right]\quad \text{ and } \quad \boldsymbol A := P\otimes I_n\: .$$   

The value of $\underline{a}$ in Assumption~\ref{ass:row_stoch} is  $0.25$. One can compute the matrix $\boldsymbol Q$ by solving  the LMI in \eqref{eq:LMI} with $\eta=0.5$ and imposing a diagonal structure. The  solution is 
\begin{equation}
\boldsymbol Q= \mathrm{diag}(0.186\,I_2,\: 0.214\, I_2,\: 0.055\, I_2,\: 0.03\, I_2)\:.
\label{eq:Q_1_static_simulation}
\end{equation}

The cost function of each agent $i$ is defined according to Assumption~\ref{eq:def_costFunction_Ji}, as
\begin{equation}
J^i(x^i,\boldsymbol x^{-i})=\frac{\gamma_i}{2}\lVert x^i-x^*_i \rVert^2_{Q_i}+\iota_{\mathcal{X}^i}(x^i) +\frac{1}{2}\lVert x^i-\boldsymbol{A} \boldsymbol{x} \rVert^2_{Q_i}\:,
\label{eq:cost_fnc_eg_1}
\end{equation}
where $Q_i\in\mathbb{R}^{2\times 2}$ is the block present in the diagonal of $\boldsymbol Q$, $x^{*,i}$ the target position of the agent, $\gamma_i\in\mathbb{R}$ the ``discover'' parameter and $\lVert x^i-\boldsymbol{A} \boldsymbol{x} \rVert^2_{Q_i}$ the aggregative term. The parameter $\gamma_i$ defines how much a robot prefers to achieve its goal position instead of staying close to the others. 
%In fact, a high value of $\gamma_i$ would lead to a final position of robot $i$ closer to the target position and vice versa. 

From the collective dynamics in (\ref{eq:group_dynamics}), the explicit update rule is obtained applying a forward-backward splitting to the operator, obtaining
\begin{multline}\label{eq:dynamics_eg_1}
\boldsymbol x(k+1) = \mathrm{proj}_{ \boldsymbol{\mathcal{X}} }
\big[ x(k)-\epsilon \left( \boldsymbol Q \big( \gamma(\boldsymbol x(k)-\boldsymbol x^*) \big) \right. \\
\left. + D(\boldsymbol x(k)-\boldsymbol A\boldsymbol x(k)) \right) \big] \:,
\end{multline}
where $D = \mathrm{diag}\left( (1-a_{i,i})_{i=1}^{N} \right)$, $a_{i,i}$'s are the diagonal elements of $\boldsymbol{A}$, $\epsilon$ is the step size of the dynamics and $\mathrm{proj}_{ \boldsymbol{\mathcal{X}}}$ is the projection operator over the set $\boldsymbol{\mathcal{X}}$.
 
\subsubsection{Simulation results}
\label{subsubsec:simulation1_result1}
In the numerical simulation, we use the following values:
\begin{itemize}
\item initial positions $x_0^1=[ 5 \: 0]^\top$, $ x_0^2=[ 20\: 0]^\top$, $x_0^3=[50 \: 0]^\top$ and $x_0^4=[10\: 0 ]^\top$;
\item desired final positions $x^{*,1}=[ 100 \: 100]^\top$, $ x^{*,2}=[ 60\: 100]^\top$, $x^{*,3}=[0 \: 50]^\top$ and $x^{*,4}=[100\: 50 ]^\top$;
\item edge of local constraint sets, $\mathcal{X}^i$, $r=5$;
\item the same discover parameter is used for all the agents, $\gamma_i = 2.5$ for all $i\in\bN[1,N]$;
\item small step size $\epsilon = 1$.
\end{itemize}

\begin{figure}[thpb]
  \centering
  \includegraphics[width=0.8\columnwidth]{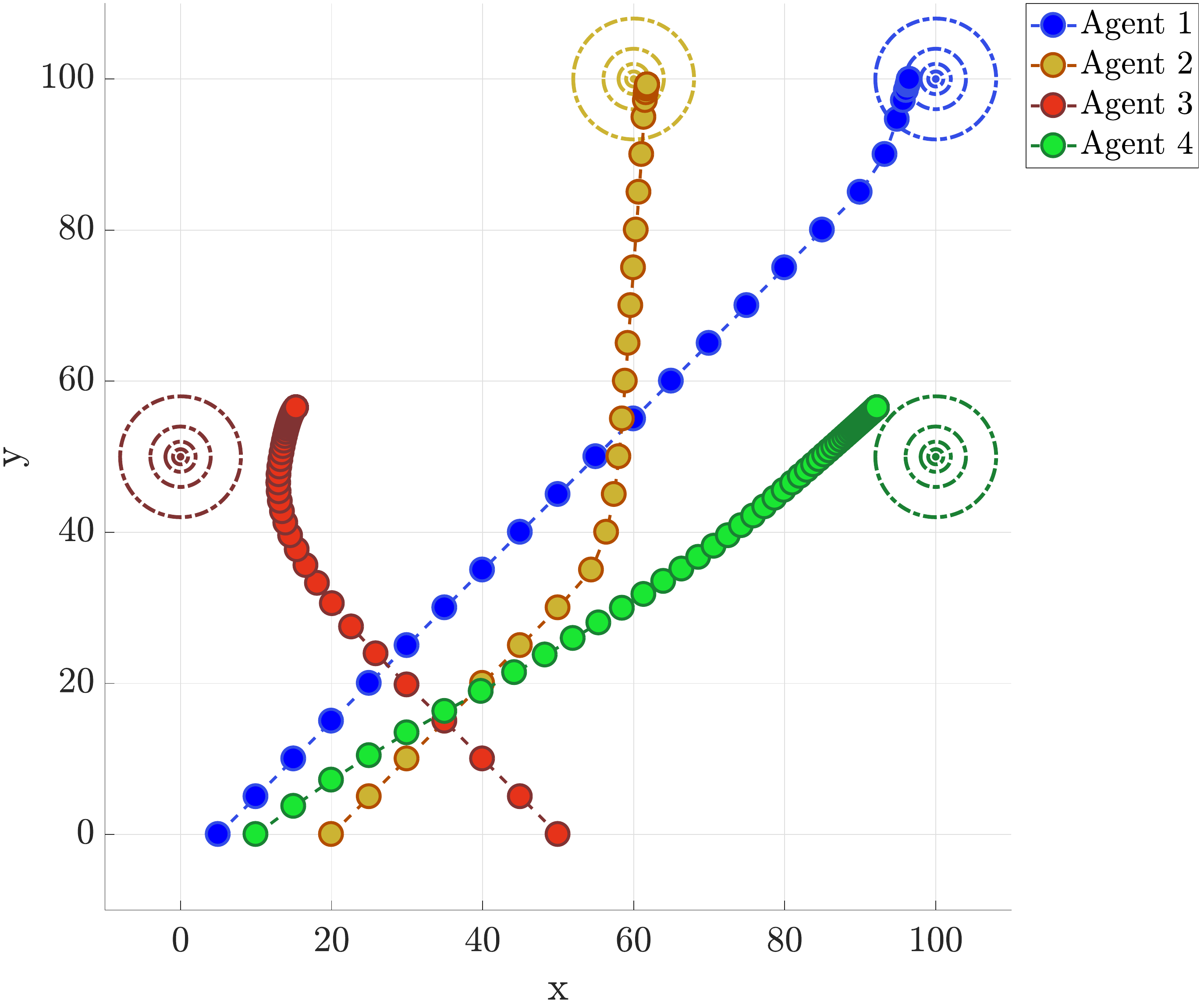}
  \caption{Trajectories of the $4$ robots generated by the dynamics in \eqref{eq:dynamics_eg_1}. The initial and desired positions are $ x_0^i$ and $x^{*,i}$, respectively. The latter represented by concentric dashed circles.}
  \label{fig:static_case}
\end{figure}

The trajectories of the four robots are shown in Fig.~\ref{fig:static_case}, where the desired position of each robot $i\in\bN[1,N]$ is represented by dashed circles with the same color of the robot $i$. 

A closer look at the matrix $P$ clarifies the behavior of each robot. The first two robots weight mostly their relative positions and, since their targets are relatively close each other, they converge very close to the desired location. The remaining two robots must instead adapt their motion also with respect to that of robots $1$ and $2$, in order to reduce the effect of the proximity term in their cost function. This leads the final positions of these robots distant from their targets.

\subsubsection{Obstacle Avoidance}
\label{subsubsec:obstacle_avoidance}
Next, we deal with an obstacle avoidance problem, which can be handled via the same setup, with the only difference that $\mathcal{X}^i$ shall be modified. When the agent $i$ approaches the object, its local constraints set is formulated such that the future robot position cannot collide with the obstacle. More precisely, we define $\mathcal A$ as the set of points covered by the obstacle, and $\hat{\mathcal{X}}^i$, the local constraint set that the agent would have without the obstacle. Finally, we define the set $\mathcal{X}^i$ as the largest convex subset of $\hat{\mathcal{X}}^i\setminus \mathcal{A}$.    

\begin{figure}[tb]
  \centering
  \includegraphics[width=0.8\columnwidth]{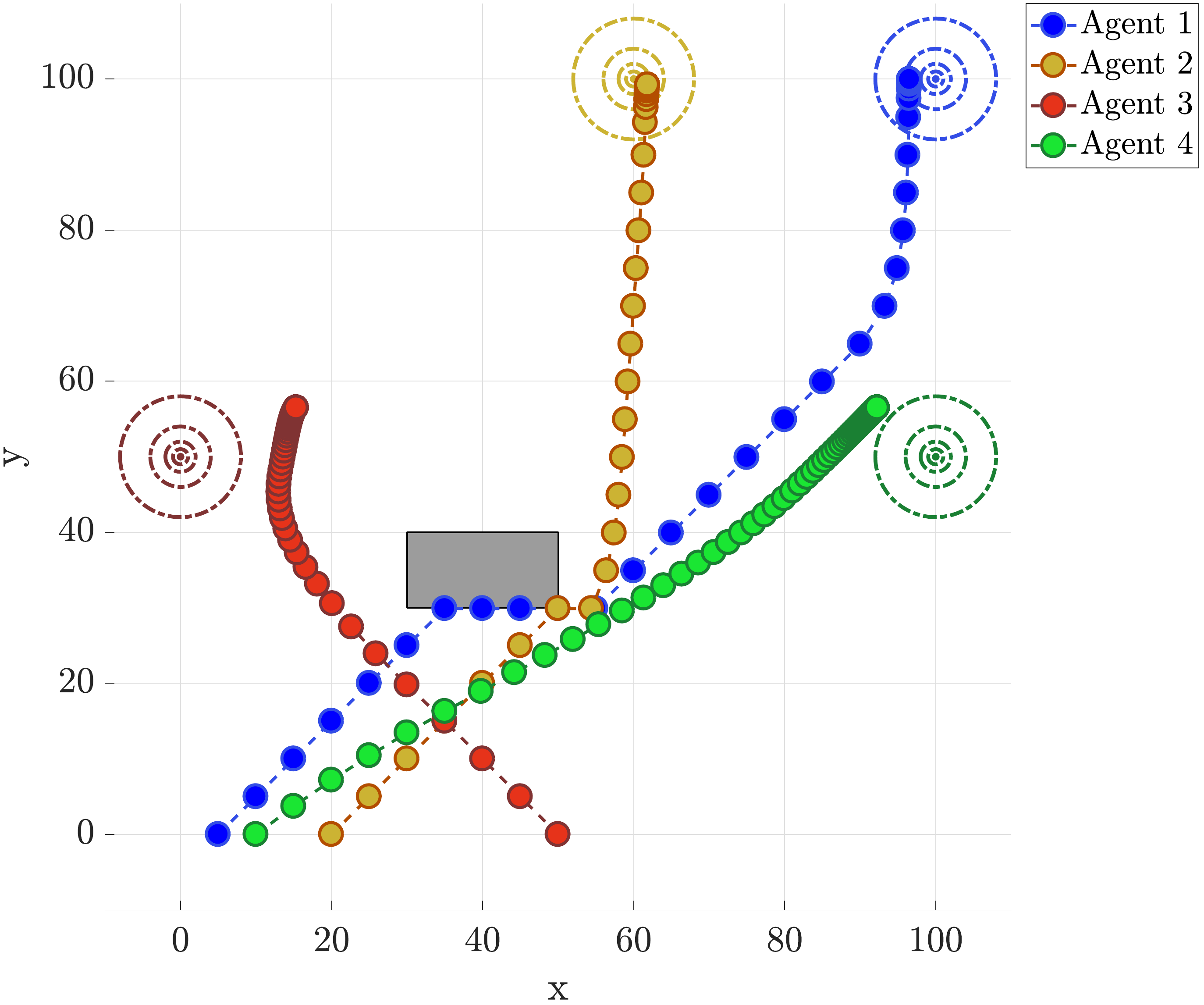}
  \caption{Trajectories of the $4$ robots generated by the dynamics in \eqref{eq:dynamics_eg_1}, performing obstacle avoidance. The obstacle is represented by the gray rectangle. The desired positions are represented by concentric dashed circles.}
  \label{fig:obtacle_avoidance}
\end{figure}

In Fig.~\ref{fig:obtacle_avoidance}, we show a resulting trajectory, run with the same parameters of the previous simulation. Each agent successfully avoids the obstacle, and the final positions almost coincide with the ones of the previous simulation.  
%Notice that the  result of the algorithm is  related to the obstacle considered and that a satisfactory result is not guaranteed. Sometimes this issue can be overcome with an \textit{ad hoc} tuning of the discovery parameters $\gamma_i$.
In our experience, expecially in the scenarios with obstacles, the tuning of the $\gamma_i$s is crucial to obtain satisfactory trajectories.

\section{Conclusion and Outlook}
\label{section:conclusions}
In this paper, we have studied the problem of a group of robots performing a distributed exploration task. We have shown that under weak condition on the communication, the global convergence to an NWE or to a PNWE can still be reached, respectively in the static and the time-varying case, even though in this latter case the imposition of a dwell time is necessary. Moreover we presented a practical implementation of the algorithm and studied its performances in different setups.
This can be seen as a work towards asynchronous proximal dynamics, which aims at highlighting the potentials and the possible applications of this research topic. 
%Multi-agent time-invariant proximal dynamics converge to a network equilibrium if the adjacency matrix of the graph is an averaged linear operator in a space with block-diagonal weight matrix. For time-varying proximal dynamics, convergence holds if the adjacency matrix switches within a set of averaged linear operators, provided that the switch is subject to a certain dwell time. Multi-agent proximal dynamics are applicable to distributed robotic applications.
Future research will investigate milder conditions and different algorithms under which proximal dynamics are guaranteed to be fully distributed. % Conditions that ensure global convergence to an equilibrium shall be explored for time-varying, switched, proximal dynamics as well.
\appendix

\subsection*{Proof of Lemma \ref{th:convergence_banach_dynamics} }

From Assumption \ref{ass:row_stoch}, $A$ is marginally stable, with no eigenvalues on the boundary of the unit disk but semi-simple eigenvalues at $1$. From  \cite[Lem. 4]{belgioioso:fabiani:blanchini:grammatico:18css}, the linear operator $A$ is averaged in some Hilbert space with norm $\tilde{Q} \succ 0$. From \cite[Def. 4.33]{Bauschke2010:ConvexOptimization} and \cite[Lem. 4]{belgioioso:fabiani:blanchini:grammatico:18css}, $\boldsymbol{A}$ is averaged in the Hilbert space with norm $\tilde{ \boldsymbol{Q} } \succ 0$. 
Then, the collective proximal operator $\boldsymbol{\mathrm{prox}}_{\boldsymbol{f}}^{\boldsymbol{Q}}(z):=\mathcal{J}_{\boldsymbol{Q^{-1}\partial f}}(z)$ is averaged in the same Hilbert space \cite[Prop. 23.34(i)]{Bauschke2010:ConvexOptimization}.
%%By Lemma \ref{lem:proximal_averaged}, the proximal operator $\boldsymbol{\textrm{prox}}^{ \boldsymbol{Q} }$ is averaged in the same Hilbert space. 
Thus, by \cite[Prop. 4.44]{Bauschke2010:ConvexOptimization}, the composition $\boldsymbol{\textrm{prox}}_{\boldsymbol{f}}^{ \boldsymbol{Q} } \circ \boldsymbol{A} $ is also averaged in $\mathcal{H}_{ \tilde{ \boldsymbol{Q} } }$. The proof then follows from \cite[Prop. 5.15]{Bauschke2010:ConvexOptimization}.
\hfill $\blacksquare$

\subsection*{Proof of Proposition \ref{prop:eta_avg_LMI_solvable} }
Follows directly from \cite[Def. 4.33]{Bauschke2010:ConvexOptimization}, \cite[Lem. 4]{belgioioso:fabiani:blanchini:grammatico:18css} and  Lemma~\ref{th:convergence_banach_dynamics}.
\hfill $\blacksquare$

\subsection*{Proof of Theorem \ref{th:dwell_time_theorem}}
For each  $i \in \bN[1,M]$, define $T_i := \boldsymbol{\mathrm{prox}}_{\boldsymbol{f}}^{\boldsymbol Q_i} \circ \boldsymbol A_i$. From  \cite[Prop. 4.44]{Bauschke2010:ConvexOptimization} we know that $T_i$, is averaged as well, say $\eta_i-$averaged. Thus, without any switching, the sequence $(\boldsymbol x(k))_{k=0}^\infty$ generated by the Banach--Picard iteration would converge to some vector in $\fix(T_i)$. 

Due to Assumptions ~\ref{ass:linear_regular}, \ref{ass:infinite_occurence_of_P} and \cite[Lemma~3.8, Fact~5.3(i)]{Bausche:convergence_regular_operators}, we have that, for all $i\in\bN [1,N]$, 
%%and $\eta_i \in (0,1-\underline a_i)$,
\begin{align}
&\mathrm{d}_{\boldsymbol Q_i}\big(\boldsymbol{x}(k+n),\mathcal{E}\big) \leq 2 \phi^k_i \mathrm{d}_{\boldsymbol Q_i}\big(\boldsymbol{x}(n),\mathcal{E}\big)
\label{eq:norm_inequality_reg_lin}
\end{align}
where we defined the parameters
\begin{align}
&\hspace{5mm}\phi_i:=\sqrt{\dfrac{\alpha_i^{-1}\kappa_i^2}{1+\alpha^{-1}_i\kappa_i^2}}\in[0,1), \ \alpha_i:=\frac{1-\eta_i}{\eta_i}. \nonumber
\end{align}

%Then the lower bound of the dwell time in \eqref{eq:claim_tau} can be computed via arguments similar to those in \cite[Sec.~3.2.1]{Liberzon2003:SwitchingInSystemAndControl}.

%%%%%%%%%%%%%%%%%%%%%%%%%%%%%%%
Then the lower bound of the dwell time in \eqref{eq:claim_tau} can be computed via arguments similar to \cite[Sec.~3.2.1]{Liberzon2003:SwitchingInSystemAndControl}.
Consider $M=2$, since the distance functions are norms, we have that
\begin{equation}
\begin{split}
{\scriptstyle \lambda_{\min,1}\mathrm{d}\big(\boldsymbol{x}(k),\mathcal{E}\big)\leq \mathrm{d}_{\boldsymbol Q_1}\big(\boldsymbol{x}(k),\mathcal{E}\big) \leq \lambda_{\max,1} \mathrm{d}\big(\boldsymbol{x}(k),\mathcal{E}\big)}\\
{\scriptstyle \lambda_{\min,2}\mathrm{d}\big(\boldsymbol{x}(k),\mathcal{E}\big)\leq \mathrm{d}_{\boldsymbol Q_2}\big(\boldsymbol{x}(k),\mathcal{E}\big) \leq \lambda_{\max,2} \mathrm{d}\big(\boldsymbol{x}(k),\mathcal{E}\big) }
\end{split}\:.
\label{eq:bound_distances}
\end{equation}
%where $\lambda_{\min,i}:=\mathrm{min}\big(\Lambda(\boldsymbol Q_i)\big)$ and $\lambda_{\max,i}:=\mathrm{max}\big(\Lambda(\boldsymbol Q_i)\big)$ denote minimum-norm and maximum-norm eigenvalues, respectively.

Now, suppose that $\sigma(k)=1$ for all $k \in \bN[k_0,k_1-1]$, and $\sigma(k)=2$ for all $k \in \bN[k_1,k_2-1]$, where $k_{j+1}-k_{j}>\tau$. 
By \eqref{eq:norm_inequality_reg_lin} and \eqref{eq:bound_distances}, we obtain
\begin{equation}
\mathrm{d}_{\boldsymbol Q_2}\big(\boldsymbol{x}(k_1),\mathcal{E}\big)\leq  2{ \frac{\lambda_{\max,2}}{\lambda_{\min,1}}}\phi_1^\tau\mathrm{d}_{\boldsymbol Q_1}\big(\boldsymbol{x}(k_0),\mathcal{E}\big),
\label{eq:dwell_calc_first_step}
\end{equation}
and consequently
\begin{equation}
\mathrm{d}_{\boldsymbol Q_1}\big(\boldsymbol{x}(k_2),\mathcal{E}\big)\leq  4{\frac{\lambda_{\max,2}\lambda_{\max,1}}{\lambda_{\min,1}\lambda_{\min,2}}}(\phi_2\phi_1)^\tau\mathrm{d}_{\boldsymbol Q_1}\big(\boldsymbol{x}(k_0),\mathcal{E}\big).
\label{eq:dwell_calc_second_step}
\end{equation}
By \cite[Th.~3.1]{Liberzon2003:SwitchingInSystemAndControl} and some manipulation, we conclude the converge to a PNWE if the dwell time $\tau$ satisfies
%\begin{align*}
%4{\scriptstyle \frac{\lambda_{\max,2}\lambda_{\max,1}}{\lambda_{\min,1}\lambda_{\min,2}}}(\phi_2\phi_1)^\tau \leq 1 \,,
%\end{align*}
%and, equivalently, if
\begin{align*}
\tau\geq \mathrm{log}_{\phi_2\phi_1}\left(\frac{\lambda_{\min,1}\lambda_{\min,2}}{4\lambda_{\max,2}\lambda_{\max,1}}\right).
\end{align*}

The proof for the general case when $M \ge 2$ is analogous and leads to the lower bound on the dwell time in \eqref{eq:claim_tau}.
%%%%%%%%%%%%%%%%%%%%%%%%%%%%%%%
\hfill $\blacksquare$

%% BIBLIOGRAPHY %%%%%%%%%%%%%%%%%%%%%%%%%%%%%%%%%%%%%%%%%%%%%%%%%%%%%%%%%%%%%
\bibliographystyle{IEEEtran}
\bibliography{diary_bibliography.bbl,librarySG.bbl}

\end{document}